\documentclass[10pt,leqno,oneside]{extarticle}
\usepackage{cite}

\usepackage{amsmath}
\usepackage{amssymb}
\usepackage{amsfonts}
\usepackage{setspace}

\usepackage{latexsym}
\usepackage{mathrsfs}
\usepackage{epsfig}

\newcommand{\n}{\noindent}
\newcommand{\be}{\begin{equation}}
\newcommand{\ee}{\end{equation}}
\newcommand{\ben}{\begin{displaymath}}

\newcommand{\een}{\end{displaymath}}
\newcommand{\ep}{\hspace{\stretch{1}}$\Box$}

\newcommand{\vs}{\vspace{0.2cm}}

\newtheorem{Remark}{Remark}

\newtheorem{Proposition}{Proposition}
\newtheorem{Theorem}{Theorem}

\newtheorem{Lemma}{Lemma}

\begin{document}

\begin{center}
{\huge A note on scalar curvature and the 

\vs
convexity of boundaries.}

\vspace{0.6cm}
{\large Martin Reiris.}\footnote{e-mail: martin@aei.mpg.de}\\

\vs
\textsc{\small Max Planck Institute f\"ur Gravitationsphysik \\ (Albert Einstein Institut)\\ Golm - Germany.}\\

\vspace{0.7cm}
\begin{minipage}[c]{10cm}
\begin{spacing}{1}
{\small We prove that any smooth Riemannian manifold of non-negative scalar curvature and with a strictly mean convex and compact boundary component can be ($C^{2}$) extended beyond the component to have non-negative scalar curvature and to enjoy anyone of the following three types of (new) boundary: strictly convex, totally geodesic or strictly concave. The extension procedure can be applied for instance to ``positive mass" type of theorems.}
\end{spacing}
\end{minipage}
\end{center}

\section{Introduction.}


\vs
We prove an isometric extension procedure for manifolds with non-negative scalar curvature and strictly mean convex, compact boundary. The extended manifolds have non-negative scalar curvature and can be made to have anyone of the following three types of boundary: strictly convex, totally geodesic or strictly concave (we recall the precise definitions below). Theorems of the sort, namely stating that every element of a class ${\mathcal A}$ of Riemannian manifolds with boundary can be isometrically embedded into an element of a special (and hopefully interesting) class ${\mathcal B}$ of Riemannian manifolds, are usually appreciated for the sole reason that they permit to import geometric properties enjoyed by elements in the class ${\mathcal B}$ to those in the class ${\mathcal A}$. An example of a procedure of the type, indeed not entirely foreign to the one discussed here, is given in \cite{MR2379775}. There, an isometric extension procedure was introduced for manifolds with boundary that preserve lower (sectional) curvature bounds and produces a totally geodesic boundary.  The tool proved to be useful to study a series of questions on convergence and compactness of Riemannian manifolds with boundary. 
%
%
We will point out some applications of the present work after recalling basic terminology and after stating the extension Theorem (Theorem \ref{MT}).  The proof is given in Section \ref{Proof}. 

A manifold $M$ is $C^{k+1},\ k\geq 0$ if the transition functions of coordinate charts are $C^{k+1}$.  A Riemannian manifold $(M,g)$ is $C^{k}$ if $M$ is $C^{k+1}$ and the metric components of $g$ in any coordinate chart are $C^{k}$. $(M,g)$ is smooth if $k=\infty$. 
We will say that a $C^{k}$ ($k\geq 0$) Riemannian manifold $(\bar{M},\bar{g})$ is an {\it extension} of a smooth Riemannian manifold $(M,g)$ if there is a $C^{k+1}$ embedding from $(M,g)$ into $(\bar{M},\bar{g})$, which is also a $C^{k}$ isometry. 

Let $(M,g)$ be a smooth Riemannian manifold. We will denote compact boundary components of $M$ by $\partial^{c}M$.  Let $\Omega\subset M$ be a region with smooth boundary $\partial \Omega$. Denote by $h$ the metric on $\partial \Omega$ inherited from $g$, by $\Theta$ the second fundamental form of $\partial \Omega$ with respect to the outgoing normal\footnote{If $\varsigma$ is the outgoing unit normal then for any two vectors $v_{1}$ and $v_{2}$ in $T_{q}\partial \Omega$ ($q\in \partial \Omega$) we have $\Theta(v_{1},v_{2})=<\nabla_{v_{1}}\varsigma,v_{2}>$.} 
 and by $\theta=tr_{h}\Theta$ the mean curvature ($tr_{h}$ is the trace with respect to $h$). Under this notation recall that $\partial^{c} M$ is strictly mean-convex (resp. mean-concave) if $\theta>0\ ({\rm resp.} <0)$ and is strictly convex (resp. strictly concave) if $\Theta>0\ ({\rm resp}. <0)$ (as a symmetric two-form, i.e. $\Theta(v,v)>0\ ({\rm resp.} <0)$ for any $v\neq 0$). $\partial^{c} M$ is totally geodesic if $\Theta=0$.

\begin{Theorem}\label{MT} {\rm (The extension Theorem).} Let $(M,g)$ be a smooth Riemannian manifold ($n=dim(M)\geq 3$) of non-negative scalar curvature and with a strictly mean convex and compact boundary component $\partial^{c} M$. Then, there are extensions beyond $\partial^{c}M$ to $C^{2}$ manifolds $(\bar{M},\bar{g})$, of non-negative scalar curvature $R$ and enjoying any of the following three types of boundaries: strictly convex, totally geodesic or strictly concave. Moreover, in any of the extensions, $\bar{M}\setminus \varphi(Int(M))$ is diffeomorphic to $\partial^{c} M\times [0,1]$.
\end{Theorem}  
In simple terms, the manifold $(M,g)$ can be extended to a ``collar" around $\partial^{c} M$ to have non-negative scalar curvature and a new boundary component, replacing $\partial^{c}M$, with either strictly convex, totally geodesic or strictly concave boundary (as wished).  As a matter of fact the manifolds $\bar{M}$ will be $C^{\infty}$, but $\bar{g}$ just $C^{2}$. 
\begin{Remark}\label{REM1} As a byproduct of the construction, $(M,g)$ can be extended beyond $\partial^{c} M$ by adding a collar, in such a way that the new boundary component $\partial^{c} \bar{M}$ has second fundamental form $\bar{\Theta}$ enjoying the lower bound
\be\label{LBU}
\bar{\Theta}>\frac{\theta_{0}}{2(n-1)}\bar{h}
\ee
where $\theta_{0}>0$ is a lower bound for the mean curvature of $\partial^{c} M$ and $\bar{h}$ is the metric on $\partial^{c} \bar{M}$ induced from $\bar{g}$. Moreover the distance from $\partial^{c} \bar{M}$ to $\partial^{c} M$ can be made as small as wished still preserving (\ref{LBU}).
\end{Remark} 
The manifold $M=[0,1]\times T^{2}$ $(T^{2}=S^{1}\times S^{1}$) with the flat metric $g=dx^{2}+d\theta_{1}^{2}+d\theta_{2}^{2}$ ($\theta_{1}$ and $\theta_{2}$ are the angular coordinates of the $S^{1}$ factors), cannot be extended beyond its boundary to have strictly convex boundary and non-negative scalar curvature. This shows that the hypothesis of strict convexity of the boundary component of $M$ cannot be removed. 

A simple example of a manifold $M$ with $R\geq 0$, strictly mean convex boundary which is neither convex, totally geodesic or concave, together with extensions $\bar{M}_{i},i=1,2,3$ with strictly convex, totally geodesic or strictly concave boundary (respectively to $i=1,2,3$) is  given in the following. Consider the unit two-sphere $S^{2}$ in polar coordinates (from a point) $\{r,\varphi\}$ and the manifold $S^{1}$ with coordinate $\{\theta\}$. Then, on $S^{2}\times S^{1}$ consider the metric $g=dr^{2}+\sin^{2} r d\varphi^{2}+(1+\epsilon\cos 4r)^{2} d\theta^{2}$. If $\epsilon>0$ is small enough $g$ has positive scalar curvature. The manifold $M=\{(r,\varphi,\theta)/0\leq r\leq \pi/6\}$ endowed with $g$ has (if $\epsilon$ is small enough) mean  convex boundary but not strictly convex boundary. The manifolds $\bar{M}_{i}=\{(r,\varphi,\theta)/0\leq r\leq r_{i}\}$ ($i=1,2,3$) with $r_{1}=\pi/3$, $r_{2}=\pi/2$ and $r_{3}=2\pi/3$, and endowed with $g$ are extensions with strictly convex, totally geodesic and strictly concave boundary respectively.   

There are a number of deformation techniques that one can find in the literature \cite{MR2844438},\cite{MR1982695} related to ``positive mass" type of theorems which share some elements with the extension theorem developed here but seem to be of a different nature. As in these works one can apply also the extension Theorem \ref{MT} to obtain ``positive mass" type of theorems.   
To illustrate this we show here how the (Riemannian) positive mass theorem for asymptotically flat manifolds with non-negative scalar curvature and with strictly mean convex compact boundaries (see \cite{MR1626060} and references therein) can be easily reduced to the standard positive mass theorem for boundary-less asymptotically flat manifold. Consider $(M,g)$ a smooth manifold of non-negative scalar curvature, asymptotically flat ends, and strictly mean convex (in the outward direction) boundary with possibly many connected components. Theorem \ref{MT} allows to extend $(M,g)$ to a manifold $(\bar{M},\bar{g})$ with totally geodesic boundary and non-negative scalar curvature. The manifold $(\bar{M},\bar{g})$ can then be ``doubled" along its boundary to obtain a boundary-less, asymptotically flat $C^{2}$ manifold (the regularity $C^{2}$ is seen easily from the construction of $\bar{M}$) $(\bar{\bar{M}},\bar{\bar{g}})$ of non-negative scalar curvature. If  $(\bar{\bar{M}},\bar{\bar{g}})$ is known to have positive mass at any one of its ends (for instance if $\bar{\bar{M}}$ is spin) then, obviously, $(M,g)$ will also have positive mass at any one of its ends. Further applications will be discussed elsewhere.

\section{Proof of the main result.}\label{Proof}

From now on we will assume, without loss of generality and to simplify notation, that $(M,g)$ is a smooth Riemannian manifold and that $\partial M$ has only one connected component and is compact. Therefore we will write $\partial M$ instead of $\partial^{c}M$. 

Below we will describe a simple and concrete setup to prove Theorem \ref{MT}. We collect it in the statement of Lemma \ref{MT2} from which Theorem \ref{MT} directly follows. A proof of the Lemma is given afterwards.

Consider the tubular neighborhoods $T(\partial M,\Gamma)$ of 
$\partial M$
\ben
T(\partial M,\Gamma):=\{p\in M/dist(p,\partial M)\leq \Gamma\}
\een
\n If $\Gamma^{-}>0$ is small enough, then for any point $p\in T(\partial M,\Gamma^{-})$ there is a unique point $q(p)\in \partial M$ such that $dist(p,\partial M)=dist(p,q(p))$ and, moreover, the map $\psi:T(\partial M,\Gamma^{-})\rightarrow [-\Gamma^{-},0]\times \partial M$ given by $\psi(p)=(-dist(p,\partial M),q(p))$ is a smooth diffeomorphism. Denote by $s(p)=-dist(p,\partial M)$ the coordinate on $[-\Gamma^{-},0]$. Also denote by $\bar{h}(s)$ the metric induced from $\psi_{*} g$ into $\{s\}\times \partial M$. Note, of course, that $\{s\}\times \partial M$ is canonically diffeomorphic to $\partial M$ and also note that for this reason $\bar{h}(s)$, $s\in [-\Gamma^{-},0]$ can be thought (as we will do) as a smooth path of metrics over $\partial M$. 
Abusing slightly notation we will write $g=ds^{2}+\bar{h}$ (this is justified because (of course)
$
(g-ds^{2})|_{(s,q)}(v_{1}+a_{1}\partial_{s},v_{2}+a_{2}\partial_{s})=\bar{h}(s)|_{q}(v_{1},v_{2})
$
%
for any $q\in \partial M$, $v_{1},v_{2}$ in $T_{q} \partial M$ and real numbers $a_{1},a_{2}$).
Finally, we will use the notation 
\begin{gather}\label{hiDef}
h_{0}:=\bar{h}\bigg|_{s=0},\ 
h'_{0}:=\frac{\partial}{\partial s} \bar{h}\bigg|_{s=0^{-}},\
h''_{0}:=\frac{\partial^{2}}{\partial s^{2}} \bar{h}\bigg|_{s=0^{-}}
\end{gather}
for $\bar{h}(s)$ and its first and second normal (one sided) derivatives at $s=0$.

Now, if we can extend the path $\bar{h}(s)$, which is so far defined in the interval $[-\Gamma^{-},0]$, to a path of metrics, (also denoted by) $\bar{h}(s)$, on the interval $[\Gamma^{-},\Gamma^{+}]$, $\Gamma^{+}>0$, in such a way that when we consider $\bar{h}(s)$ restricted only to $[0,\Gamma^{+}]$, it is a $C^{2}$ path of metrics satisfying: 

\vs
I. the Riemannian metric $ds^{2} +\bar{h}$ (over $[0,\Gamma^{+}]\times \partial M$) has non-negative scalar curvature,

\vs
II. $
h_{0}=\bar{h}\bigg|_{s=0},\ 
h'_{0}=\frac{\partial}{\partial s} \bar{h}\bigg|_{s=0^{+}},\ 
h''_{0}=\frac{\partial^{2}}{\partial s^{2}} \bar{h}\bigg|_{s=0^{+}}$,

\vs  
III. $\frac{\partial}{\partial s} \bar{h}\bigg|_{s=\Gamma^{+}}>0\ ({\rm resp.}\ =0,<0)$,

\vspace{0.3cm}
\n then, the metric $\bar{g}=ds^{2}+\bar{h}$, over $[-\Gamma^{-},\Gamma^{+}]\times \partial M$, will be $C^{2}$, will have non-negative scalar curvature, and moreover, as
\ben
\Theta\bigg|_{s=\Gamma^{+}}=\frac{1}{2}\frac{\partial}{\partial s} \bar{h}\bigg|_{s=\Gamma^{+}}
\een
then (from III) the boundary component $\{\Gamma^{+}\}\times \partial M$ will be strictly convex (resp. totally geodesic or strictly concave). In such case the Riemannian manifold  
\ben
(M,g)\cup_{\psi}([-\Gamma^{-},\Gamma^{+}]\times \partial M,\bar{g})
\een

\n where $\cup_{\psi}$ identifies, via $\psi$, the region $T(\partial M,\Gamma^{-})$ on $M$ with the region $[-\Gamma^{-},0]\times \partial M$ on $[-\Gamma^{-},\Gamma^{+}]\times \partial M$, would be the desired extension claimed in the Theorem \ref{MT}. It follows that Theorem \ref{MT} is a direct consequence of the following Lemma.

\begin{Lemma}\label{MT2} Let $(M,g)$ be a smooth Riemannian manifold of non-negative scalar curvature and mean convex, compact and connected boundary. Let $h_{0},h'_{0},h''_{0}$ be as in (\ref{hiDef}). Then there is a $C^{2}$ path of metrics $\bar{h}(s),s\in [0,\Gamma^{+}]$, for some $\Gamma^{+}>0$, over $\partial M$ and satisfying {\rm I-II-III}.
\end{Lemma}

\vs
\n {\it Proof:} We will concentrate to obtain an extension $(\bar{M},\bar{g})$ with strictly convex boundary (case $>0$ in III) and we will indicate at the end how to obtain the other two possibilities.

We will use $h_{0}$ as a background metric to estimate expressions derived from symmetric two-tensors. In particular we will use the following (usual) inner product. Let $U$ be a symmetric two-tensor field on $\partial M$. Given $q\in \partial M$ let $\{e_{i}(q)\}$ be any $h_{0}$-orthonormal basis of tangent vectors. Let $U_{ij}:=U(e_{i},e_{j})$. Then $<U,V>_{0}=\sum_{i,j} U_{ij}V_{ij}$, is a pointwise ($\{e_{i}\}$-invariant) inner product on symmetric two-tensor fields. Write $|U|^{2}_{0}(q)=<U(q),U(q)>_{0}$. For any $U$ and tangent vectors $v,w$ (at the same point $q$) we have $|U(v,w)|\leq |U|_{0}|v|_{0}|w|_{0}$ where $|v|_{0}$ ($|w|_{0}$) denotes the $h_{0}$-norm of $v$ ($w$). For any Riemannian metric $h$ on $\partial M$ let $tr_{h} U$ denote the trace of $U$ with respect to $h$. We have $tr_{h} U=U_{ij} h^{ij}$ where $h^{ij}$ is the inverse matrix to $h_{ij}$. Defining $h^{-1}$ by $h^{-1}_{ij}=h^{ij}$ we have $tr_{h}U=<U,h^{-1}>_{0}$. Also we will denote by $|U|_{h}$ the (pointwise) norm of $U$ but with respect to $h$.

It will be more convenient to define $\bar{h}(s)$ in terms of a $C^{2}$ path of metrics $h(t)$, where $s$ and $t$ are related by
$ds=\alpha(t)dt$, and where $\alpha(t)$ is a positive $C^{1}$ function of $t$, required to satisfy $\alpha(0)=1$ and $(d\alpha/dt)(0)=0$ (these two conditions are important) that will be chosen conveniently later. In other words given $h(t)$ and $\alpha(t)$ define $\bar{h}(s)$ by $\bar{h}(s):=h(t(s))$, where $t(s)$ would be found by inverting $s(t)=\int_{0}^{t} \alpha(\bar{t})d\bar{t}$ (but we won't need to do so). Thus
\be\label{gmDef}
ds^{2}+\bar{h}=\alpha^{2}dt^{2}+h,
\ee   
Note that $s(t)$ is a $C^{2}$ function of $t$ (and $t(s)$ a $C^{2}$ function of $s$) if $\alpha(t)$ is $C^{1}$ and therefore $ds^{2}+\bar{h}(s)=ds^{2}+h(t(s))$ is $C^{2}$ (over any chart $(s,x_{1},\ldots,x_{n-1})$ with $(x_{1},\ldots,x_{n-1})$ a chart on $\partial M$, $n=\dim(M)$). The remark is important as at the end the function $\alpha(t)$ we will be $C^{1}$ but not $C^{2}$.  

The idea now is to define $h(t)$ independently of $\alpha(t)$ to satisfy automatically the following two conditions:

\vs
II'. $h_{0}=h\bigg|_{t=0},\ h'_{0}=\frac{\partial}{\partial t} h\bigg|_{t=0^{+}},\ 
h''_{0}=\frac{\partial^{2}}{\partial t^{2}} h\bigg|_{t=0^{+}}$,

\vs  
III'. $\frac{\partial}{\partial t} h\bigg|_{t=t^{+}}>0\ (=0,<0)$,

\vs
\n and then note that if II' and III' hold so do II and III (for $\bar{h}(s)$) independently of the function $\alpha(t)$ defining $t(s)$. Then for the given $h(t)$, find a function $\alpha(t)$ to satisfy I. In this form $\bar{h}(s)$ will satisfy I, II and III. As a matter of fact $h(t)$ and $\alpha(t)$ will be given at the end explicitly in terms of $h_{0},h'_{0},h''_{0}$ and some constants defined out of them. More precisely $h(t)$ will be defined in (\ref{hDef}) and $\alpha(t)$ will be defined as (\ref{Aa}) (for some value of $a$ explained later) over an interval $[0,t^{I}]$ and as (\ref{Ab}) (for some value of $b$ explained later) over $[t^{I},t^{+}]$.  

Derivatives with respect to $s$ will be denoted with a tilde ($'$) and with respect to $t$ with a dot ($\dot{}$).  

We define $h(t)$ by
\be\label{hDef}
h(t)=h_{0}+\frac{(tr_{h_{0}}h'_{0})t}{n-1} h_{0} + F(t) \hat{h'_{0}} + G(t) h''_{0}
\ee 

\begin{figure}[h]
\centering
\includegraphics[width=8cm,height=8cm]{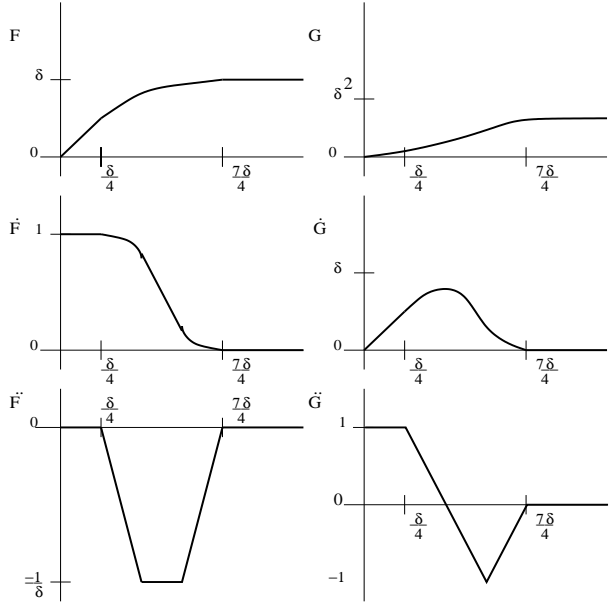}
\caption{Sketches of $F$ and $G$ and their first and second derivatives.}
\label{SFG}
\end{figure}

\n where the hat $\hat{}$ in $h'_{0}$ denotes the traceless part of $h'_{0}$ (with respect to $h_{0}$) and $F(t)$ and $G(t)$ are two $C^{2}$ functions dependent on a parameter $\delta>0$ (fixed later) which are precisely described in what follows. The function $F(t)$ is defined 
to be $C^{2}$, to have $F(0)=0$, $\dot{F}(0)=1$ and with second derivative given by
\be\label{FSD}
\ddot{F}(t)=\left\{ 
\begin{array}{lcl}
0 & \text{if} & t\in [0,\frac{\delta}{4}]\cup [\frac{7\delta}{4},\infty),\vspace{0.15cm}\\
-\frac{2}{\delta^{2}}(t-\frac{\delta}{4}) & \text{if} & t\in [\frac{\delta}{4},\frac{3\delta}{4}],\vspace{0.15cm}\\
-\frac{1}{\delta} & \text{if} & t\in [\frac{3\delta}{4},\frac{5\delta}{4}],\vspace{0.15cm}\\
-\frac{1}{\delta}+\frac{2}{\delta^{2}}(t-\frac{5\delta}{4}) & \text{if} & t\in [\frac{5\delta}{4},\frac{7\delta}{4}]
\end{array}
\right.
\ee
The explicit form of $F(t)$ can be found by integrating (\ref{FSD}) twice, but it is of no importance here. The function $G(t)$ is defined to be $C^{2}$, to have $G(0)=0$, $\dot{G}(0)=0$ and with second derivative given by 
\be\label{GSD}
\ddot{G}(t)=\left\{
\begin{array}{lcl}
1  & \text{if} & t\in [0,\frac{\delta}{4}],\vspace{0.15cm}\\
1-\frac{2}{\delta}(t-\frac{\delta}{4}) & \text{if} & t\in [\frac{\delta}{4},\frac{5\delta}{4}],\vspace{0.15cm}\\
-1+\frac{2}{\delta}(t-\frac{5\delta}{4})&\text{if} & t\in [\frac{5\delta}{4},\frac{7\delta}{4}],\vspace{0.15cm}\\
0 & \text{if} & t\in [\frac{7\delta}{4},\infty)
\end{array}
\right.
\ee
Again, the explicit form of $G(t)$ can be found by integrating twice (\ref{GSD}) but it is of no importance here. A sketch of $F(t)$ and $G(t)$ can be seen in Figure \ref{SFG}. Most of what quantitatively matters to us are the following simple global bounds (for all $t$)
\begin{gather*}
|F|\leq \delta,\ |\dot{F}|\leq 1,\ |\ddot{F}|\leq \frac{1}{\delta},\\
|G|\leq \delta^{2},\ |\dot{G}|\leq \delta,\ |\ddot{G}|\leq 1
\end{gather*}
along with the explicit expressions (which can be easily deduced)
\begin{gather*}
F(t)=t,\ \dot{F}(t)=1,\ \ddot{F}(t)=0,\\
G(t)=\frac{t^{2}}{2},\ \dot{G}(t)=t,\ \ddot{G}(t)=1
\end{gather*}
for $t$ in $[0,\frac{\delta}{4}]$, and 
\begin{gather*}
F(t)=\delta,\ \dot{F}(t)=0,\ \ddot{F}(t)=0,\\
G(t)=\frac{1}{2}\frac{47}{48}\delta^{2},\ \dot{G}(t)=0,\ \ddot{G}(t)=0
\end{gather*}
for $t$ in $[\frac{7\delta}{4},\infty)$.

\vs
Observe that $h_{0}+U$, with $U$ a symmetric two-tensor field, is point-wise positive definite, namely a Riemannian metric, as long as (pointwise) $|U|_{0}\leq 1/2$ (if $|v|_{0}=1$ then $h_{0}(v,v)+U(v,v)\geq 1-|U|_{0}\geq 1/2>0$). Therefore from (\ref{hDef}), the bounds $|F|\leq \delta$ and $|G|\leq \delta^{2}$, and the mean convexity hypothesis $2\theta_{0}=tr_{h_{0}}h{'}_{0}>0$, we obtain that $h(t)$, for all $t\geq 0$, will be a $C^{2}$ path of metrics provided $\delta\leq \delta_{0}$ with $\delta_{0}$ such that
\ben
\delta_{0}|\hat{h}'_{0}|_{0}+\delta_{0}^{2}|h''_{0}|_{0}\leq \frac{1}{2},\ {\rm (pointwise)}
\een 
{\it From now on we will assume $\delta_{0}$ was fixed and that $\delta\leq \delta_{0}\leq 1$}. On the other hand, as $\dot{F}(t)=0$ and $\dot{G}(t)=0$ for $t\geq 7\delta/4$, then, for $t\geq 7\delta/4$ we have
\ben
\dot{h}(t)=\frac{tr_{h_{0}}h'_{0}}{n-1}h_{0}
\een  
Therefore, the second fundamental forms $\Theta$ of the slices $\{t\}\times \partial M$, for any $t\geq 7\delta/4$, which are given by $\Theta=\dot{h}/2\alpha$, are positive definite regardless of the (positive) values of $\alpha(t)$. Summarizing: we define $h(t)$ by (\ref{hDef}) which with $t^{+}\geq 7\delta/4$ satisfies II' and III' automatically. 
%
Lemma \ref{MT2} will be thus proved as long as we can chose $\delta$ ($\delta\leq \delta_{0}$) and find $\alpha(t)$ defined at least over an interval $[0,t^{+}]$, with $t^{+}\geq 7\delta/4$, for which the metric $\alpha^{2}dt^{2}+h$ has non-negative scalar curvature $R$, namely I holds. We pass now to explain how to find such $\delta$ and $\alpha(t)$.

\begin{figure}[h]
\centering
\includegraphics[width=6cm,height=6cm]{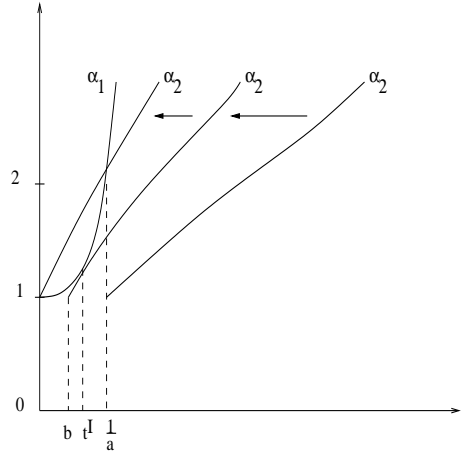}
\caption{Sketches of $\alpha_{1}$, $\alpha_{2}$ and the intermediate time $t^{I}$. There are three graphs of $\alpha_{2}$ corresponding to three different values of the parameter $b$: $b=1/a$ (right), final $b$ (middle), $b=0$ (left). As $b\downarrow 0$ the graphs of the functions $\alpha_{2}$ shift to the left, this is indicated by the arrows.}
\label{A}
\end{figure}

First we observe that the scalar curvature $R$ of (\ref{gmDef}) over every slice $\{t\}\times \partial M$ can be written as
\be\label{SCE}
\alpha^{2}R=\frac{\dot{\alpha}}{\alpha} (tr_{h}\dot{h}) - (tr_{h}\ddot{h}) + \frac{3}{4} |\dot{h}|^{2}_{h} - \frac{1}{4} (tr_{h}\dot{h})^{2} + \alpha^{2} {\mathcal{R}}
\ee      
where ${\mathcal{R}}$ is the scalar curvature of $h(t)$ as a metric over $\{t\}\times \partial M\ (\sim \partial M)$. To see this, derivate $2\alpha \theta= tr_{h}\dot{h}=\dot{h}_{ij}h^{ij}$ with respect to $t$ and then use the expressions (valid in any dimension greater or equal than $3$)\footnote{To obtain the first derivate $h_{ik}h^{kj}=\delta_{i}^{\ j}$, the second is the well known Riccati equation and to obtain the third contract twice the {\it Gauss}-equation (\!\!\cite{MR2371700}, pg. 38).}  
\begin{gather*}
(h{^{ij}})\dot{}=-(\dot{h}_{lm})h{^{li}}h{^{mj}},\\
\dot{\theta}=-\alpha(|\Theta|^{2}_{h}+Ric(n,n)),\\
Ric(n,n)+|\Theta|_{h}^{2}= \frac{1}{2}(R+|\Theta|_{h}^{2}+\theta^{2}-{\mathcal{R}}),\\
|\Theta|^{2}_{h}=\frac{|\dot{h}|^{2}_{h}}{4\alpha^{2}}
\end{gather*} 
Now, $R$ will be non-negative on a domain $[0,t^{+}]\times \partial M$ if for every $t\in [0,t^{+}]$ the right hand side of (\ref{SCE}), as a function on $\partial M$, is non-negative. We will think this condition on the non-negativity of the right hand side of (\ref{SCE}) as a condition on $\alpha(t)$. The function $\alpha(t)$ that will satisfy this condition will be defined separately as a function (that we call) $\alpha_{1}$ on the interval $[0,t^{I}]$ and as a function (that we call) $\alpha_{2}$ on $[t^{I},t^{+}]$ that will be seen to match $C^{1}$ at $t=t^{I}$ ($t^{I}$ accounts for ``intermediate" time). Let us give a glimpse of what will come to better orient the reading. First, as can be seen from the expressions (\ref{Ab}) and (\ref{Aa}), the functions $\alpha_{1}$ and $\alpha_{2}$ will depend (after $\delta,c_{1},c_{2},c_{3}$ have been fixed, see later) on parameters $a$ and $b$ respectively. Moreover it will be that if $a\geq a_{0}\geq 4/\delta$ and $0<b<\delta/4$ the expression (\ref{SCE}) with $\alpha=\alpha_{1}$ or $\alpha=\alpha_{2}$ will be non-negative for $t$ on the intervals $[1,1/a]$ and $[b,\Gamma^{+}]$ respectively.  Now, the required values of $a$ and $b$ to fix $\alpha_{1}$ and $\alpha_{2}$ and therefore $\alpha$, will be such that the graphs of (\ref{Ab}) and (\ref{Aa}) touch tangentially at a point $\alpha_{1}(t^{I})=\alpha_{2}(t^{I})$ for some $t^{I}$ such that $b<t^{I}<1/a$. A representation is given in Figure (\ref{A}). We will be explaining all this is more in detail in what follows.

At the moment we move to construct $\alpha_{2}$ and to justify its properties. We will use the explicit expression (\ref{hDef}) that we have for $h(t)$ to get a suitable lower bound expression for the right hand side of (\ref{SCE}) in terms of $\alpha$, $\delta$ and constants $\bar{c}_{0},\bar{c}_{1},\bar{c}_{2}$ depending only on $h_{0}, h'_{0}, h''_{0}$, and then find $\alpha_{2}(t)$, over a suitable interval $[b,t^{+}]$ to make such lower bound (with $\alpha=\alpha_{2}$) non-negative (zero indeed). The expression for the referred lower bound will be obtained from the following proposition and given explicitly afterwards in (\ref{LB}). 

\begin{Proposition}\label{P1} There are $\delta_{1}\leq \delta_{0}$, $t_{1}\leq 1$, and positive numbers $\bar{c}_{0},\ \bar{c}_{1},\ \bar{c}_{2}$ depending on $(h_{0},h'_{0},h''_{0})$ such that for any $\delta\leq \delta_{1}$ and $0\leq t\leq t_{1}$ we have the pointwise bounds (over $\partial M$)
\begin{gather}
\label {E1} tr_{h}\dot{h}\geq \frac{1}{2} (tr_{h_{0}}h'_{0})\geq \bar{c}_{0},\\
\label{E2} |tr_{h}\ddot{h}-\frac{3}{4} |\dot{h}|^{2}_{h}+\frac{1}{4} (tr_{h}\dot{h})^{2}|\leq \frac{\bar{c}_{1}}{\delta},\\
\label{E3} {\mathcal{R}(h)}\geq -\bar{c}_{2}
\end{gather}

\end{Proposition}

\n {\it Proof of the Proposition \ref{P1}:}  We show first (\ref{E1}). We need a couple of observations. First, from the expression (\ref{hDef}) and the bounds $|F|\leq \delta$, $|G|\leq \delta^{2}$, we observe that for any $\epsilon>0$ there are $\bar{\delta}_{1}\leq \delta_{0}$ and $\bar{t}_{1}\leq 1$ such that for any $0\leq \delta\leq \bar{\delta}_{1}$ and $0\leq t\leq \bar{t}_{1}$ we have (pointwise) $|h-h_{0}|_{0}\leq \epsilon$, (here and below we make of course $h=h(\delta,t)$). Second, it is simple to see that there is $\epsilon_{0}$ such that for any $\epsilon\leq \epsilon_{0}$, if $|h-h_{0}|_{0}\leq \epsilon$ then $|h^{-1}-h_{0}|_{0}\leq 2\epsilon$. From this and the general inequalities for symmetric two-tensor fields $U$ 
\ben
|tr_{h_{0}}U|-|h^{-1}-h_{0}|_{0}|U|_{0}\leq |tr_{h}U|\leq |tr_{h_{0}}U|+|h^{-1}-h_{0}|_{0}|U|_{0}
\een
we deduce that if $|h-h_{0}|\leq \epsilon\leq \epsilon_{0}$ then for any $U$ we have
\be\label{trl}
|tr_{h_{0}}U|-2\epsilon |U|_{0}\leq |tr_{h}U|\leq |tr_{h_{0}}U|+2\epsilon |U|_{0}
\ee
Combing the two observations we have obtained that for any $\epsilon\leq \epsilon_{0}$ there are $\bar{\delta}_{1}(\epsilon)\leq \delta_{0}$ and $\bar{t}_{1}(\epsilon)\leq 1$ such that for any $0\leq \delta\leq \bar{\delta}_{1}(\epsilon)$, $0\leq t\leq \bar{t}_{1}(\epsilon)$ and symmetric two-tensor field $U$, the inequalities (\ref{trl}) hold.

Now taking the trace of the time derivative of (\ref{hDef}) we get
\ben
tr_{h}\dot{h}=\frac{1}{n-1}(tr_{h_{0}}h'_{0})(tr_{h}h_{0})+\dot{F}(tr_{h}\hat{h}'_{0})+\dot{G} (tr_{h}h''_{0})
\een
Using (\ref{trl}) with $U=h_{0},\hat{h}'_{0},h''_{0}$ we conclude that for any $\epsilon\leq \epsilon_{0}$ there are $\bar{\delta}_{1}(\epsilon)\leq \delta_{0}$ and $\bar{t}_{1}(\epsilon)\leq 1$ such that for any $0\leq \delta\leq \bar{\delta}_{1}(\epsilon)$, $0\leq t\leq \bar{t}_{1}(\epsilon)$ we have
\ben
tr_{h}\dot{h}\geq (tr_{h_{0}}h'_{0}) - \frac{2\epsilon}{\sqrt{n-1}}(tr_{h_{0}}h'_{0}) - 2\epsilon |\hat{h}'_{0}|_{0} - \delta |tr_{h_{0}}h''_{0}|_{0} -2\epsilon\delta |h''_{0}|_{0}
\een 
In particular if $\delta\leq \delta_{1}\leq \bar{\delta}_{1}(\epsilon)$ and $t\leq t_{1}\leq \bar{t}_{1}(\epsilon)$ then 
\be\label{FE}
tr_{h}\dot{h}\geq (tr_{h_{0}}h'_{0}) - \frac{2\epsilon}{\sqrt{n-1}} (tr_{h_{0}}h'_{0}) - 2\epsilon |\hat{h}'_{0}|_{0} - \delta_{1} |tr_{h_{0}}h''_{0}|_{0} -2\epsilon \delta_{1} |h''_{0}|_{0}
\ee 
Choosing now $\epsilon(\leq \epsilon_{0})$ sufficiently small and then $\delta_{1}(\leq \bar{\delta}_{1}(\epsilon))$ and $t_{1}(\leq \bar{t}_{1}(\epsilon))$ sufficiently small we deduce from (\ref{FE}) that for any $0\leq \delta \leq \delta_{1}$ and $0\leq t\leq t_{1}$
\ben
tr_{h}\dot{h}\geq \frac{1}{2} (tr_{h_{0}}h'_{0})\geq \inf \frac{1}{2}(tr_{h_{0}}h'_{0}):=\bar{c}_{0}>0.
\een
as desired (the right hand side defines $\bar{c}_{0}$).
  
To obtain (\ref{E2}) we proceed similarly. Taking the trace of the second time derivative of (\ref{hDef}) and using the bounds $|\ddot{F}|\leq 1/\delta$, $|\ddot{G}|\leq 1$ we obtain
\ben
|tr_{h}\ddot{h}|=|\ddot{F}(tr_{h}\hat{h}'_{0})+\ddot{G}(tr_{h}h''_{0})|\leq \frac{|tr_{h}\hat{h}_{0}'|}{\delta}+|tr_{h}h''_{0}|
\een  
We use the same $\delta_{1}$ and $t_{1}$ as was chosen for (\ref{E1}) before. Then using (\ref{trl}) with $U=\hat{h}'_{0},h''_{0}$ in the previous equation we obtain for any $0\leq \delta\leq \delta_{1}$ and $0\leq t\leq t_{1}$
\ben
|tr_{h}\ddot{h}|\leq \frac{\bar{c}_{3}}{\delta} + \bar{c}_{4}
\een
where $\bar{c}_{3}$, $\bar{c}_{4}$ depend only on $(h_{0},h'_{0},h''_{0})$. Along similar lines one obtains $-(3/4)|\dot{h}|^{2}_{h}+(1/4)(tr_{h}\dot{h})^{2})\leq \bar{c}_{5}$ depending only on the data $(h_{0},h'_{0},h''_{0})$. Therefore
\ben
tr_{h}\ddot{h}-\frac{3}{4} |\dot{h}|^{2}_{h}+\frac{1}{4} (tr_{h}\dot{h})^{2}\leq \frac{\bar{c}_{3}+\delta(\bar{c}_{4}+\bar{c}_{5})}{\delta}\leq \frac{\bar{c}_{3}+\bar{c}_{4}+\bar{c}_{5}}{\delta}:=\frac{\bar{c}_{1}}{\delta}
\een
as desired (the last equality defines $\bar{c}_{1}$ and we used $\delta\leq \delta_{0}\leq 1$).

Finally  (\ref{E3}) follows from the fact that the biparametric family of $C^{2}$ metrics $h(\delta,t)$, $0\leq \delta \leq \delta_{0}, 0\leq t \leq 1$ (on $\partial M$) given in (\ref{hDef}), depends continuously (in $C^{2}$) on $(\delta,t)$ and that the set $\{(\delta,t)/ 0\leq \delta \leq \delta_{0}, 0\leq t \leq 1\}$ is compact. 
Therefore for any $0\leq \delta\leq \delta_{0}$ and $0\leq t\leq 1$ and consequently for any $0\leq \delta\leq \delta_{1}$ and $0\leq t\leq t_{1}$ ($\delta_{1}$ and $t_{1}$ as chosen for (\ref{E1}) and (\ref{E2}) before) we have ${\mathcal{R}}\geq -\bar{c}_{5}$ depending only on $(h_{0},h'_{0},h''_{0})$. \ep

\vs
The estimations (\ref{E1})-(\ref{E3}) lead, for any $\alpha>0$ with $\dot{\alpha}>0$ to the following lower bound to the right hand side of (\ref{SCE})
\be\label{LB}
\frac{\dot{\alpha}}{\alpha}\bar{c}_{0} - \frac{\bar{c}_{1}}{\delta}-\bar{c}_{2}\alpha^{2}
\ee
as long as $\delta\leq \delta_{1}$ and $t\leq t_{1}$. Thus, $R\geq 0$ if $\alpha>0$, $\dot{\alpha}>0$ and $\dot{\alpha}/\alpha - c_{1}/\delta -c_{2}\alpha^{2}\geq 0$ where we have defined $c_{1}=\bar{c}_{1}/\bar{c}_{0}$ and $c_{2}=\bar{c}_{2}/\bar{c}_{0}$. 

\vs
$\bullet$ The function $\alpha_{2}$ is defined by
\be\label{Aa}
\alpha_{2}(t)=\sqrt{\frac{c_{1}}{\delta}\frac{1}{\big( (\frac{c_{1}}{\delta}+c_{2})e^{-2(t-b)\frac{c_{1}}{\delta}}-c_{2}\big)}}
\ee
It depends on a parameter $b$ that we require to lie in $(0,\delta/4)$ and that will be fixed later. If $\delta_{1}$ and $t_{1}$ are chosen sufficiently small then we claim that (\ref{Aa}) is well defined (i.e. there are no zeros in the denominator) at least on the interval $[b,t^{+}:=4\delta]$ ({\it note that we are defining $t^{+}:=4\delta$}). Indeed, if $b\leq t\leq \Gamma^{+}=4\delta$ we have the estimate
\ben
(\frac{c_{1}}{\delta}+c_{2})e^{-2(t-b)\frac{c_{1}}{\delta}}-c_{2}\geq \frac{c_{1}}{\delta} e^{-8c_{1}}-c_{2}
\een
where $\delta\leq \delta_{1}$. If $\delta_{1}$ is sufficiently small then the right hand side of the previous expression is strictly positive.  Now we further impose $0<\delta\leq \inf\{t_{1},\delta_{1}\}/4$. With this condition on $\delta$ (and the definition of $\Gamma^{+}$ as $4\delta$) (\ref{LB}) is a lower bound for the right hand side of (\ref{SCE}) for any function $\alpha>0$ with $\dot{\alpha}>0$ defined on $[b,\Gamma^{+}]$. But $\alpha_{2}$ verifies 
\begin{gather}
\alpha_{2}(b)=1,\\
\label{AaI} \dot{\alpha_{2}}/\alpha_{2} - c_{1}/\delta -c_{2}\alpha_{2}^{2}= 0
\end{gather}
from which it follows that $\alpha_{2}(t)>0$ and $\dot{\alpha_{2}}(t)>0$ on $[b,t^{+}]$ and that with the choice $\alpha=\alpha_{2}$ we have $R\geq 0$ over $[b,t^{+}]\times \partial M$, no matter the value of $b$ in  $(0,\delta/4)$. 

\vs
$\bullet$ The function $\alpha_{1}$ is defined by 
\be\label{Ab}
\alpha_{1}(t)=1+a^{2}t^{2}
\ee
It depends on a parameter $a$ that we require to lie in $(4/\delta,\infty)$ and that will be fixed later. Note that $1\leq \alpha_{1}\leq 2$ for $t$ in $[0,1/a]$. We claim that there is $a_{0}\geq 4/\delta$ such that if $a\geq a_{0}$ (but no matter which value), the right hand side of (\ref{SCE}) with $\alpha=\alpha_{1}$ will be non-negative and thus $R\geq 0$ over $[0,1/a]\times \partial M$. Indeed, one has expansions 
\begin{gather*}
tr_{h}\dot{h}=tr_{h_{0}}h'_{0}+O(t),\\
-tr_{h}\ddot{h}+ \frac{3}{4} |\dot{h}|^{2}_{h} - \frac{1}{4} (tr_{h}\dot{h})^{2}=-tr_{h_{0}}\ddot{h_{0}}+ \frac{3}{4} |\dot{h_{0}}|^{2}_{h_{0}} - \frac{1}{4} (tr_{h_{0}}\dot{h_{0}})^{2} +O(t),\\  
{\mathcal{R}}={\mathcal{R}}(h_{0}) + O(t)
\end{gather*}
where each of the $O(t)$ is of the form $O(t)=t f + O(t^{2})$ with $f$ a function on $\partial M$ depending only on $(h_{0},h'_{0},h''_{0})$ and $|O(t^{2})|\leq c t^{2}$ with $c$ a constant also dependent on $(h_{0},h'_{0},h''_{0})$.
Evaluating (\ref{SCE}) at $t=0$ we get that the scalar curvature $R$ at the initial time $\{0\}\times \partial M$ (which we will denote by $R_{0}$ below) is given by $R_{0}=-tr_{h_{0}}\ddot{h_{0}}+ \frac{3}{4} |\dot{h_{0}}|^{2}_{h_{0}} - \frac{1}{4} (tr_{h_{0}}\dot{h_{0}})^{2}+{\mathcal{R}}(h_{0})$. With this information we can write the right hand side of (\ref{SCE}) as
\be\label{EFIN}
\frac{\dot{\alpha}}{\alpha}(tr_{h_{0}}h'_{0}+O(t)) + R_{0} +O(t) \alpha^{2} + (\alpha^{2}-1){\mathcal{R}}_{0}+O(t)
\ee  
with $O(t)$ as explained before. Now, by the strict mean convexity hypothesis we have $tr_{h_{0}}h'_{0}>0$ and also by the non-negative scalar curvature hypothesis we have $R_{0}\geq 0$. Finally $\dot{\alpha}_{1}=2a^{2} t$ and $\alpha^{2}_{1}-1=(2+a^{2}t^{2})a^{2}t^{2}$. Making $z=at$ then $\dot{\alpha}_{1}=2z/a$ and $\alpha^{2}_{1} -1=(2+z^{2})z^{2}$. It is clear then that if 
$a\geq a_{0}$ with $a_{0}$ big enough (depending on $(h_{0},h'_{0},h''_{0})$, we require also $a_{0}\geq 4/\delta$) then the expression (\ref{EFIN}) and therefore the right hand side of (\ref{SCE}) will be (pointwise) non-negative for $t$ in the interval $[0,1/a]$ (i.e. $z\in [0,1]$). This shows the claim.

We proceed now to adjust $a$ and $b$ to fix $\alpha_{1}$ and $\alpha_{2}$ and therefore fix $\alpha$. First we make an observation. 
Consider the family of functions $\{\alpha_{2}\}$ parametrized by $b\in [0,\delta/4]$. Recall that for each $b$, the corresponding function $\alpha_{2}$ is defined on $[b,\Gamma^{+}]$. It is simple to see from (\ref{Aa}) and (\ref{AaI}) that there is a uniform bound (i.e. independent of $b\in [0,\delta/4]$) for $|\alpha_{2}|$, $|\dot{\alpha}_{2}|$. Moreover, for all $b\in [0,\delta/4]$ the corresponding function $\alpha_{2}$ has positive derivative ($\dot{\alpha}_{2}>0$) on its domain $[b,\Gamma^{+}]$, in particular for $b=0$. 
Now fix $a\geq a_{0}$ in such a way that $\dot{\alpha}_{1}(t=1/a)=2a$ is strictly greater than the uniform bound for $|\dot{\alpha}_{2}|$. With this we will fix $b$ and $t^{I}$ as follows. Note first that the graph of $\alpha_{1}$ on $[0,1/a]$ and the graph of $\alpha_{2}$ for $b=1/a$ on $[b=1/a,\Gamma^{+}]$ obviously do not intersect.  Now, starting from $b=1/a$, decrease $b$ and consider the graphs of the functions $\alpha_{2}(t),\ t\in [b,\Gamma^{+}]$ which shift to the left as $b\downarrow 0$. Then, because of how $a$ was chosen, we have that for some $0<b<1/a$ there is $t^{I}\in (b,1/a)$ for which 
\begin{gather*}
\alpha_{1}(t^{I})=\alpha_{2}(t^{I}),\\
\dot{\alpha}_{1}(t^{I})=\dot{\alpha}_{2}(t^{I})
\end{gather*}
A picture of this can be seen in Figure \ref{A}. This fixes $b$ and $t^{I}$. Thus we have defined $\alpha_{1}$, $\alpha_{2}$ and $t^{I}$ and therefore the function $\alpha(t)$. This finishes the construction and the proof of Lemma \ref{MT2} for strictly convex boundaries. 

\vs

We show now how to obtain the other two boundary possibilities, namely totally geodesic and strictly concave. We discuss first the totally geodesic case. We will use part of the construction of the previous strictly convex case. In particular take the same $\delta$ and $\alpha$ as we did before but change in the expression (\ref{hDef}) for $h(t)$ the factor $t$ (of $tr_{h_{0}}h'_{0})h_{0}$) for example by the $C^{2}$ function $H(t)$
\ben
H(t)=\left\{
\begin{array}{lcl}
t & \text{ for }&  t\in [0,3\delta],\vs \\
\frac{\ddot{H}_{\epsilon}}{6\epsilon}(t-\delta)^{3}+(t-3\delta)+3\delta & \text {for } & t\in [3\delta,3\delta+\epsilon],\vs \\
\epsilon^{1/4}\dot{H}_{\epsilon}\sin\frac{(t-3\delta-\epsilon)}{\epsilon^{1/4}}+H_{\epsilon} \cos \frac{(t-3\delta-\epsilon)}{\epsilon^{1/4}} & \text{ for } & t\in [3\delta+\epsilon,t^{+}]
\end{array}
\right.
\een  
where $H_{\epsilon}=H(3\delta+\epsilon)$, $\dot{H}_{\epsilon}=\dot{H}(3\delta+\epsilon)$, $\ddot{H}_{\epsilon}=\ddot{H}(3\delta+\epsilon)=-(1/\epsilon^{1/2})H_{\epsilon}$ are given by
\begin{gather*}
H_{\epsilon}=\frac{3\delta+\epsilon}{1+\frac{\epsilon^{3/2}}{6}},\
\dot{H}_{\epsilon}=1-\frac{\epsilon^{1/2}}{2}(\frac{3\delta+\epsilon}{1+\frac{\epsilon^{3/2}}{6}}),\
\ddot{H}_{\epsilon}=-\frac{1}{\epsilon^{1/2}}(\frac{3\delta+\epsilon}{1+\frac{\epsilon^{3/2}}{6}})
\end{gather*}
where $t^{+}=3\delta+\epsilon+\epsilon^{1/4}\arctan (\dot{H}_{\epsilon}\epsilon^{1/4}/H_{\epsilon})$. Note that now $t^{+}$ changed and is no more equal to $4\delta$. The small number $\epsilon>0$ will be adjusted below. That the boundary $\{t=t^{+}\}$ is totally geodesic can be seen directly from the fact that $\dot{h}(t^{+})=(tr_{h_{0}}h'_{0})\dot{H}(t^{+})h_{0}$ and that $\dot{H}(t^{+})=0$. That the scalar curvature of (\ref{gmDef}) is non-negative can be seen as follows. First on $[0,3\delta]\times \partial M$ the metric is the same that we have constructed in the strictly convex case which had non-negative scalar curvature and indeed positive at $\{t=3\delta\}$. 
Observe for this that the right hand side of (\ref{SCE}) is, for $t\in [2\delta,4\delta]$, greater or equal than $\frac{\dot{\alpha}}{\alpha}\bar{c}_{0}+(\frac{\dot{\alpha}}{\alpha}\bar{c}_{0} - \frac{\bar{c}_{1}}{\delta}-\bar{c}_{2}\alpha^{2})$ which because of (\ref{LB}) is greater or equal than $\frac{\dot{\alpha}}{\alpha}\bar{c}_{0}>0$. Thus the scalar curvature $R$ is positive on $[2\delta,4\delta]\times \partial M$.
Second, on $[3\delta,3\delta+\epsilon]\times\partial M$, $\ddot{H}\leq 0$ so the second term on the right hand side of (\ref{SCE}) is positive. Therefore from (\ref{SCE}) we have
\be\label{SCEII}
\alpha^{2} R\geq \frac{\dot{\alpha}}{\alpha} (tr_{h}\dot{h}) + \frac{3}{4} |\dot{h}|^{2}_{h} - \frac{1}{4} (tr_{h}\dot{h})^{2} + \alpha^{2} {\mathcal{R}}
\ee
Moreover, for $t$ in $[3\delta,3\delta+\epsilon]$ it is 
\ben
1-\frac{\epsilon^{1/2}}{2}(\frac{3\delta+\epsilon}{1+\frac{\epsilon^{3/2}}{6}})\leq \dot{H}(t)\leq 1
\een
and consequently the range of $\dot{H}(t)$ (as $t$ varies over $[3\delta,3\delta+\epsilon]$) converges uniformly to one as $\epsilon$ tends to zero. Thus as $\epsilon$ approaches zero the right hand side of (\ref{SCEII}) converges (pointwise) to
the scalar curvature $R$ at $\{t=3\delta\}\subset [0,3\delta]\times \partial M$, which as we pointed out before, is positive. It follows that for $\epsilon$ small enough the scalar curvature $R$ over $[3\delta,3\delta+\epsilon]\times \partial M$ is positive. Third, on $[3\delta+\epsilon,t^{+}]\times \partial M$ we have $\ddot{H}=-(1/\epsilon^{1/2})H$ and moreover, in this range of $t$
\ben
H_{\epsilon}\leq H(t)\leq \epsilon^{1/4}\dot{H}_{\epsilon}+H_{\epsilon},\ 0\leq \dot{H}(t)\leq 1-\frac{\epsilon^{1/2}}{2}(\frac{3\delta+\epsilon}{1+\frac{\epsilon^{3/2}}{6}})
\een
Therefore, if $\epsilon$ is small enough and recalling that $\ddot{h}=(tr_{h_{0}}\dot{h}_{0})\ddot{H} h_{0}$, the second term on the right hand side of (\ref{SCE}), which is positive, dominates all the other terms thus giving non-negative scalar curvature (positive actually).

The extension with strictly concave boundary is obtained simply by considering the same $h(t)$ and $\alpha(t)$ as we constructed above for the totally geodesic case but on the slightly bigger interval $[0,t^{+}+\bar{\epsilon}]$ with $\bar{\epsilon}$ is a sufficiently small number.
\ep


\bibliographystyle{plain}

\bibliography{Master}

\begin{thebibliography}{1}

\bibitem{MR2371700}
Arthur~L. Besse.
\newblock {\em Einstein manifolds}.
\newblock Classics in Mathematics. Springer-Verlag, Berlin, 2008.
\newblock Reprint of the 1987 edition.

\bibitem{MR2844438}
Simon Brendle and Fernando~C. Marques.
\newblock Scalar curvature rigidity of geodesic balls in {$S^n$}.
\newblock {\em J. Differential Geom.}, 88(3):379--394, 2011.

\bibitem{MR1626060}
Marc Herzlich.
\newblock The positive mass theorem for black holes revisited.
\newblock {\em J. Geom. Phys.}, 26(1-2):97--111, 1998.

\bibitem{MR1982695}
Pengzi Miao.
\newblock Positive mass theorem on manifolds admitting corners along a
  hypersurface.
\newblock {\em Adv. Theor. Math. Phys.}, 6(6):1163--1182 (2003), 2002.

\bibitem{MR2379775}
Jeremy Wong.
\newblock An extension procedure for manifolds with boundary.
\newblock {\em Pacific J. Math.}, 235(1):173--199, 2008.

\end{thebibliography}

\end{document}